{\rm}\input amstex
\documentstyle{amsppt}
\NoBlackBoxes
\NoRunningHeads
\magnification=\magstep1
\pageheight{6.8in}

\def\bs{\bigskip}
\
\def\fp{\flushpar}
\def\un{\underbar}

\TagsOnRight

\topmatter
\title A four parameter generalization of G\"ollnitz's (big) 
partition theorem\endtitle
\author {K. Alladi, G. E. Andrews, and A. Berkovich}\footnote {Research 
of the first and third authors was supported in part by a grant from the 
Number Theory Foundation and of the second author by a grant from the 
National Science Foundation}\endauthor
\affil 
{\it{Dedicated to our friend Barry McCoy on his sixtieth birthday}}
\endaffil
\subjclass Primary 05A15,05A19,11P81, 11P83\endsubjclass
\keywords partitions, G\"ollnitz theorem, four parameter key identity, 
q-series 
\endkeywords
\abstract
We announce a new four parameter partition theorem from which the (big) 
theorem of G\"ollnitz follows by setting any one of the parameters equal 
to 0. This settles a problem of Andrews who asked whether there exists 
a result that goes beyond the partition theorem of G\"ollnitz. We state 
a four parameter q-series identity (key identity) which is the generating 
function form of this theorem. In a subsequent paper, the proof of the 
new four parameter key identity will be given.
\endabstract 
\endtopmatter

\centerline{\bf{\S1. Introduction}}
\bs

Our purpose here is to announce the following new partition theorem:

\proclaim{Theorem 1} Let $P(n)$ denote the number of partitions of $n$ into 
distinct parts $\equiv-2^3,-2^2,-2^1,-2^0$ (mod 15).

Let $G(n)$ denote the number of partitions of $n$ into parts $\not\equiv 
2^0,2^1,2^2,2^3$ (mod 15) such that the difference between the 
non-multiples of 15 is $\ge15$ with equality only if a part is relatively 
prime to 15, parts which are not relatively prime to 15 are $>15$, the 
difference between the multiples of 15 is $\ge60$, and the smallest 
multiple of 15 is
$$\aligned
&\ge 30+30\tau,\text{ if 7 is a part, and }\\
&\ge 45+30\tau,\text{otherwise},\endaligned
$$
where $\tau$ is the number of non-multiples of 15 in the partition.  Then
$$
G(n)=P(n).
$$
\endproclaim

While it is obvious that this result is a partition theorem of the 
Schur-G\"ollnitz type, it is not clear that it lies beyond the (big) theorem 
of G\"ollnitz [15].  It is possible to obtain a four parameter refinement 
of Theorem 1; this is stated as Theorem 2 in \S2.  From this it follows 
that Theorem 1 (and Theorem 2) generalize the G\"ollnitz theorem in the same 
sense that G\"ollnitz's theorem extends Schur's 1926 partition theorem [16].  
Thus the question raised by Andrews [12] nearly 30 years ago whether there 
exists a partition theorem that goes beyond the (big) theorem of 
G\"ollnitz, is now answered in the affirmative.

Our purpose here is only to announce the new results and describe how they 
extend G\"ollnitz's theorem.  The proof of Theorem 2 (and consequently of 
Theorem 1) will be given in full in a subsequent paper [6].

In 1995, Alladi, Andrews, and Gordon [5] obtained a three parameter 
refinement of G\"ollnitz's theorem by the use of colored partitions.  Because 
colors were labelled by letters and the integers occurring in the colors 
were indicated by subscripts (weights), this approach was called 
{\it{the method of weighted words}}.  In \S2 we will describe an 
extension of this method that leads to the new four parameter Theorem 2.  
In doing so some essentially new ideas are required. In \S2 we will 
also state the generating function form of Theorem 2 which we call a four 
parameter {\it{key identity}}. By setting any one of the parameters equal to 0 
in Theorem 2, we get the three parameter refinement of G\"ollnitz's theorem 
due to Alladi-Andrews-Gordon (see \S3). Similarly by setting any one of the 
parameters equal to 0 in the new key identity (2.6), the three parameter key 
identity for G\"ollnitz's theorem in [5] falls out.  Finally in \S4 we 
conclude with a brief description of some problems for future research 
opened up by Theorem 2.

The two parameter key identity for Schur's theorem due to Alladi and Gordon 
[9] is essentially equivalent to the $q-$Chu-Vandermonde summation as shown 
by Alladi and Berkovich [8].  The three parameter key identity for 
G\"ollnitz's theorem due to Alladi, Andrews, and Gordon [5] that extends 
the identity in [9] is substantially deeper, and its proof utilizes either the 
${}_6\psi_6$ summation of Bailey as in [5], or Jackson's $q-$analog of 
Dougall's summation as in [4].  The proof of the new four parameter key 
identity (2.6) also relies on the ${}_6\psi_6$ summation, but requires 
several new ideas and is significantly deeper than the proofs in [4] and 
[5].  That is why the proof of the new four parameter identity will be 
presented separately [6].
\bs
\centerline{\bf{\S2. Colored reformulation and a four parameter refinement}}
\bs

We consider the integer 1 occurring in four primatry colors $\Bbb A, \Bbb 
B,\Bbb C$, and $\Bbb D$, and integers $n\ge2$ occurring in these four 
primary colors as well as in the six secondary colors $\Bbb A\Bbb B,\Bbb 
A\Bbb C,\Bbb A\Bbb D,\Bbb B\Bbb C,\Bbb B\Bbb D$, and $\Bbb C\Bbb D$.  We 
assume that integers $n\ge 4$ occur in the quaternary color $\Bbb A\Bbb 
B\Bbb C\Bbb D$ in addition to the ten colors above.  The crucial thing is 
that we discard all ternary colors $\Bbb  A \Bbb B \Bbb C, \Bbb A \Bbb B 
\Bbb D, \Bbb A \Bbb C \Bbb D$, and $ \Bbb B \Bbb C \Bbb D$.

The integer $n$ in color $ \Bbb A$ is denoted by the symbol $ \Bbb A_n$ 
with similar interpretation for $ \Bbb B_n,\dots, \Bbb C \Bbb D_n, \Bbb A 
\Bbb B \Bbb C \Bbb D_n$.  In order to discuss colored partitions, we need 
an ordering among the symbols, and so we assume that
$$
\cases \text{ if $m<n$ as ordinary (uncolored) integers,}\\
\text{then $m$ in any color $<n$ in any color, and}\\
\text{if two equal integers appear in different colors, the order is given 
by}\\
\Bbb  A \Bbb B \Bbb C \Bbb  D< \Bbb A \Bbb B< \Bbb A \Bbb C < \Bbb A \Bbb 
D< \Bbb A< \Bbb B \Bbb C< \Bbb B \Bbb D< \Bbb B< \Bbb C \Bbb D< \Bbb  C<
\Bbb D.\endcases\tag2.1
$$

Next, given any partition into parts occurring in the eleven colors above, 
we denote by $a$ the number of parts in color $\Bbb A$ (the frequency of 
$\Bbb A$), with similar interpretation for $b,c$, and $d$. Pursuing the 
same notation, $ab$ will denote the number of parts in color $\Bbb A \Bbb B$, 
with similar interpretation for $ac,\dots, cd$.  {\underbar{Note that $ab$ is 
not $a$ times $b$!}} Finally, $Q$ denotes the number of parts in color $ \Bbb A \Bbb B 
\Bbb C \Bbb D$.

We are now in a position to state our main result..

\proclaim{Theorem 2} Let $i, j,k, l$ be given nonnegative integers.

Let $P(n;i,j,k,l)$ denote the number of partitions of $n$ into parts 
occurring in the four primary colors, parts in the same color being 
distinct, and with $i$ parts in color $ \Bbb A$, $j$ parts in color $ \Bbb 
B$, $k$ parts in color $ \Bbb C$, and $l$ parts in color $ \Bbb D$.

Let $G(n;a,b,c,d,ab,\dots,cd,Q)$ denote the number of partitions of $n$ 
into colored parts occurring in the frequencies as indicated, and such that 
the difference between the nonquaternary parts is $\ge1$, with equality 
only if parts are either of the same primary color, or if the larger part 
occurs in a color of higher order as given in (2.1), and the gap between 
the quaternary parts is $\ge4$, with the added condition that the least 
quaternary part is
$$\aligned
&\ge 3+2\tau,\text{ if } \Bbb A_1\text{ is a part},\\
&\ge 4+2\tau,\text{ otherwise,}\endaligned\tag2.2
$$
where $\tau$ is the number of nonquaternary parts.  Then
$$
P(n;i,j,k,l)=\sum_{\text{constraints}} G(n;a,b,c,d,ab,\dots,cd,Q).\tag2.3
$$
where the summation is over the variables $a,b,...,cd, Q$ satisfying the 
constraints 
$$\aligned
&i=a+ab+ac+ad+Q\\
&j=b+ab+bc+bd+Q\\
&k=c+ac+bc+cd+Q\\
&l=d+ad+bd+cd+Q.\endaligned\tag2.4
$$\endproclaim

A strong four parameter refinement of Theorem 1 follows from Theorem 2 upon 
replacing
$$\cases
\Bbb A_n\mapsto 15n-8, \Bbb B_n\mapsto 15n-4, \Bbb C_n\mapsto 15n-2, \Bbb 
D_n\mapsto 15n-1,\text{ for } n\ge1,\\
\text{and consequently } \Bbb A \Bbb B_n\mapsto 15n-12, \Bbb A \Bbb 
C_n\mapsto 15n-10, \Bbb A \Bbb D_n\mapsto 15n-9,\\
\Bbb B \Bbb C_n\mapsto 15n-6, \Bbb B \Bbb D_n\mapsto 15n-5, \Bbb C \Bbb 
D_n\mapsto 15n-3,\text{ for } n\ge2,\\
\text{and } \Bbb A \Bbb B \Bbb C \Bbb D_n\mapsto 15n-15,\text{ for } n\ge 
4.\endcases \tag2.5
$$

The nice thing about the substitutions (2.5) is that the ordering (2.1) becomes
$$
7<11<13<14<18<20<21<22<24<25<26<27<28<29<33<,\dots
$$
the natural ordering among the integers $\not\equiv 2^0, 2^1, 2^2, 2^3$
(mod 15).  These substutitions imply that the primary colors correspond to 
the residue classes
$$
-2^3,-2^2,-2^1,-2^0(\text{mod}\quad 15).
$$
Since
$$
2^0+2^1+2^2+2^3=15,
$$
the ternary colors correspond to the residue classses
$$
2^0,2^1,2^2,2^3(\text{mod}\quad 15).
$$
These are the four residue classes not considered in Theorem 1.  Also, 
since the residue classes relatively prime to 15 are 
$2^0,2^1,2^2,2^3,-2^3, -2^2,-2^1,-2^0$(mod 15), it follows that the secondary 
colors correspond to the nonmultiples of 15 which are not 
relatively prime to 15.  Finally, the quaternary color corresponds to the 
multiples of 15 which are $\ge45$.  These features make Theorem 1 
particularly appealing.

The generating function form of Theorem 2 is the following remarkable four 
parameter {\it{key identity}}: If $T_n=n(n+1)/2$, and 
$\tau=a+b+c+d+ab+\dots+cd$, then
$$\aligned
&\sum_{\text{constraints}}\frac{q^{T_\tau+T_{ab}+\dots+T_{cd}-bc-bd-cd+4T_{Q-1}+3Q+2Q\tau}}
{(q)_a(q)_b(q)_c(q)_d(q)_{ab}\dots (q)_{cd}(q)_Q} \times\\
&\left\{(1-q^a)+q^{a+bc+bd+Q}(1-q^b)+q^{a+bc+bd+Q+b+cd}\right\}\\
&=\frac{q^{T_i+T_j+T_k+T_l}}{(q)_i(q)_j(q)_k(q)_l},\endaligned\tag2.6
$$
where the constraints are as in (2.4) and the summation is over 
$a,b,\dots,Q$. In (2.6) we have made use of the standard notation
$$
(A)_n=(A;q)_n=
\cases
\prod^{n-1}_{j=0}(1-Aq^j), \qquad if \qquad n>0,\\
1, \qquad \qquad \qquad \qquad if \qquad n=0,\\
\prod^{-n}_{j=1}{(1-Aq^{-j})^{-1}}, \, if \qquad n<0.
\endcases
$$

In [6] we do not prove Theorem 2 combinatorially.  Instead we prove (2.6) 
using $q-$series techniques and show that (2.6) is equivalent to Theorem 2.
\bs
\centerline{\bf{\S3. Reduction to G\"ollnitz and Schur}}
\bs

If any one of the parameters $i,j,k,l$, is set equal to 0, Theorem 2 
reduces to Theorem 2 of [5], the three parameter refinement of the 
colored version of G\"ollnitz's theorem. Note that in this case, the 
quaternary color does not occur at all. We state the result with $l=0$.

\proclaim{Theorem A}  Let $P(n;i,j,k)$ denote the number of partitions of 
$n$ into parts occuring in three primary colors $ \Bbb A, \Bbb B, \Bbb C$, 
parts in the same color being distinct, with $i$ parts in color $\Bbb A$, 
$j$ parts in color $\Bbb B$, and $k$ parts in color $ \Bbb C$.

Let $G(n;a,b,c,ab,ac,bc)$ denote the number of partitions of $n$ into parts 
occuring in colors $ \Bbb A, \Bbb B, \Bbb C, \Bbb A \Bbb B, \Bbb A\Bbb C, 
\Bbb B\Bbb C$, with indicated frequencies $a, b, ..., bc$, such that the 
difference between the parts is $\ge1$ with equality only if parts are 
either of the same primary color, or if the larger part occurs in a color 
of higher order as indicated in (2.1).  Then
$$
\sum\Sb i=a+ab+ac\\ j=b+ab+bc\\ k=c+ac+bc.\endSb G(n;a,b,c,ab,ac,bc)
=P(n;i,j,k).
$$\endproclaim

In Theorem A, if we use the substitutions
$$\cases
\Bbb A_n\mapsto 6n-4, \Bbb B_n\mapsto 6n-2, \Bbb C_n\mapsto 6n-1,\text{ for 
} n\ge1,\\
{\text {yielding}}\quad \Bbb A\Bbb B\mapsto 6n-6, \Bbb A\Bbb C_n\mapsto 6n-5, 
\Bbb B\Bbb C_n\mapsto 6n-3,\text{ for } n\ge2,\endcases\tag3.1
$$
we get a three parameter refinement of the following theorem of G\"ollnitz 
[15].

\proclaim{Theorem G}  Let $P(n)$ denote the number of partitions of 
$n$  into distinct parts $\equiv -2^2, -2^1, -2^0$ (mod $6$). 

Let $G(n)$ denote the number of partitions of $n$ into parts $\ne 1$ or $3$,  
such that the difference between the parts $\ge6$, with equality only if a 
part is $\equiv-2^2,-2^1,-2^0$ (mod $6$).  Then
$$
G(n)=P(n).
$$\endproclaim

If any one of the parameters $i,j,k,l$ is set equal to 0,  then (2.6) 
reduces to the three parameter key identity for G\"ollnitz's theorem in 
[5]. For example, with $l=0$, (2.6) reduces to
$$
\sum\Sb i=a+ab+ac\\ j=b+ab+bc\\ k=c+ac+bc\endSb
\frac{q^{T_\tau +T_{ab}+T_{ac}+T_{bc-1}}\{1-q^a+q^{a+bc}\}} 
{(q)_a(q)_b(q)_c(q)_{ab}(q)_{ac} (q)_{bc}}= 
\frac{q^{T_i+T_j+T_k}}{(q)_i(q)_j(q)_k}.\tag3.2
$$
Here $\tau=a+b+c+ab+ac+bc$. If we further set one of $i,j,k$ equal to 0, 
say $k=0$, in (3.2), then we get the two parameter key identity for 
the colored version of Schur's theorem due to Alladi and Gordon [9], 
namely,
$$
\sum\Sb i=a+ab\\ j=b+ab\endSb 
\frac{q^{T_{a+b+ab}+T_{ab}}}{(q)_a(q)_b(q)_{ab}}= 
\frac{q^{T_i+T_j}}{(q)_i(q)_j}\tag3.3
$$

In extending the refined Schur theorem in [9] to the refined G\"ollnitz 
theorem in [5], the statement of the extension was routine once the 
theorems were phrased in the language of primary and secondary colors. The 
principal reason for the increase in difficulty in going up from the Schur 
theorem to the G\"ollnitz theorem is because in the lexicographic ordering, 
one of the secondary colors $\Bbb B\Bbb C$ is of higher order than  
the primary color $\Bbb A$. In addition, even though the refined G\"ollnitz 
theorem uses three primary colors, the ternary color $\Bbb A\Bbb B\Bbb C$ 
is dropped and so only a proper subset of the complete alphabet of colors 
is used. Thus the G\"ollnitz theorem is an extension of Schur's theorem in 
a direction different from the one taken by Andrews [10], [11], who in 
retrospect used a complete alphabet of colors. In going beyond the G\"ollnitz 
theorem to Theorem 2, there is a significant increase in depth and 
complexity for a variety of reasons.  The quaternary color enters in a 
rather unusual way - the quaternaries do not directly interact with primaries 
and secondaries. The only interaction between quaternaries and the other 
colors is through the lower bound imposed on the quaternaries.  This, combined 
with the uncertainity of selecting a proper subset of colors from the 
complete alphabet of four primaries, six secondaries, four ternaries, and 
one quaternary, was perhaps the reason that the solution to the problem of 
Andrews [12] remained elusive for so long.
\bs
\centerline{\bf{\S4.  Problems for investigation}}
\bs
If the expressions in (2.6) are multiplied of $A^iB^jC^kD^l$ and summed 
over $i,j,k,l$, we get on the right hand side the quadruple infinite product
$$
\prod^\infty_{m=1} (1+Aq^m)(1+Bq^m)(1+Cq^m)(1+Dq^m),\tag4.1
$$
with four free parameters $A,B,C,D$.  This opens up several avenues of 
exploration a few of which we briefly indicate here.

Theorem 2 may be viewed as a base level undilated version of Theorem 
1.  More precisely, the generating function form of Theorem 1 may be viewed 
as emerging out of (2.6) and (4.1) under the transformation
$$
\cases
\text{(dilation) } q\mapsto q^{15},\\
\text{(translations) } A\mapsto Aq^{-8}, B\mapsto Bq^{-4}, C\mapsto Cq^{-2}, 
D\mapsto Dq^{-1}.\\
\endcases\tag4.2
$$
The size of the modulus 15 and the choice of the translations involving 
powers of 2 ensures that the  colors in Theorem 2 translate into distinct 
residue classes mod 15.  If a dilation smaller than $q\mapsto q^{15}$ is 
used,  then the residue classes would overlap, and so we would be counting 
parts with weights attached.  Alladi [1] has studied weighted partition 
identities in general and discussed certain interesting reformulations of 
G\"ollnitz's theorem and their applications [2], [3], emerging out of 
{\it{small}} dilations of (3.2). In a similar spirit it would be worthwhile 
to study weighted partition theorems emerging out of (2.6) by the 
use of dilations $q\mapsto q^M$, with $M<15$. We anticipate obtaining new and 
different versions of partition theorems that have arisen in the study of 
affine  Lie algebras, and representations of symmetric groups, by such 
weighted reformulations of Theorem 1.

Recently, Alladi and Berkovich [7] have obtained the following double 
bounded version of (3.3):
$$
\sum_{k\ge 0} q^{T_{i+j-k}+T_k}\left[\matrix M-i-j+k\\ 
k\endmatrix\right] \left[\matrix M-j\\ i-k\endmatrix\right] \left[\matrix 
L-i\\ j-k\endmatrix\right]=\left[\matrix L\\ j\endmatrix\right] 
\left[\matrix M-j\\ i\endmatrix\right]q^{T_i+T_j}.\tag4.3
$$
In (4.3), the symbols $\left[\matrix n+m\\ n\endmatrix\right]$ are the 
$q-$binomial coefficients given by
$$
\left[\matrix{n+m}\\{n}\endmatrix\right]_q =
\left[\matrix{n+m}\\{n}\endmatrix\right] =
\cases
\frac{(q^{m+1})_n}{(q)_n}, \qquad if \quad n\ge0,\\
0, \qquad \qquad \quad if \quad n<0.
\endcases
$$
If we let $L,M\to\infty$ in (4.3) and make the identifications $a=i-k$, 
$b=j-k$, and $ab=k$, then we get (3.3). 

In addition, Alladi and Berkovich have also obtained a double bounded 
version of (3.2), namely,
$$\aligned
&\sum_{constraints} q^{T_\tau+T_{ab}+T_{ac}+T_{bc-1}}\\
&\qquad \bigg\{q^{bc}
\left[\matrix L-\tau+a\\ a\endmatrix\right] 
\left[\matrix L-\tau+b\\ b\endmatrix\right] 
\left[\matrix M-\tau+c\\ c\endmatrix\right] 
\left[\matrix L-\tau\\ ab\endmatrix\right] 
\left[\matrix M-\tau\\ ac\endmatrix\right]
\left[\matrix M-\tau\\ bc\endmatrix\right]\\
&\qquad +\left[\matrix L-\tau+a-1\\ a-1\endmatrix\right]
\left[\matrix L-\tau+b\\ b\endmatrix\right] 
\left[\matrix M-\tau+c\\c\endmatrix\right] 
\left[\matrix L-\tau\\ ab\endmatrix\right] 
\left[\matrix M-\tau\\ ac\endmatrix\right] 
\left[\matrix M-\tau\\ bc-1\endmatrix\right]\bigg\}\\
&=\sum_{t\ge0} q^{t(M+2)-T_t+T_{i-t}+T_{j-t}+T_{k-t}}\left[\matrix L-t\\ 
t\endmatrix\right] \left[\matrix L-2t\\ i-t\endmatrix\right] \left[\matrix 
L-i-t\\ j-t\endmatrix\right] \left[\matrix M-i-j\\ k-t\endmatrix\right], 
\endaligned
\tag4.4
$$
where the constraints are as in (3.2) and $\tau=a+b+c+ab+ac+bc$. The proof 
(4.4), which is quite intricate, is given in [8]. If we let 
$L,M\to\infty$, then (3.2) follows because only the term corresponding to 
$t=0$ on the right hand side in (4.4) makes a contribution.  Based on the 
discovery of (4.4), we now ask whether there exists a finite bounded 
version of (2.6) (which reduces to (2.6) when certain parameters tend to 
infinity).

In the last decade, many new generalizations of the Rogers-Ramanujan 
identities were  discovered and proved by McCoy and collaborators (see [14] 
for a review and references), using the so called  thermodynamic Bethe 
ansatz (TBA) techniques.  It would be highly desirable to find a TBA 
interpretation of the new identity (2.6).  Such an interpretation, besides 
being of substantial interest in physics, may provide insight into how to 
extend Theorem 2 to five or more primary colors.
\bs
\fp\un{Acknowledgements}:  K.A. and A.B. would like to thank Carl 
Pomerance for support and encouragement, and Mel Nathanson and the other 
organizers of the DIMACS conference for the invitation to present this work.
\bs

\bigskip
\noindent
Department of Mathematics, University of Florida, Gainesville, FL 32611

\noindent
alladi\@math.ufl.edu
\bigskip
\noindent
Department of Mathematics, The Pennsylvania State University, University 
Park, PA 16802

\noindent
andrews\@math.psu.edu
\bigskip
\noindent
Department of Mathematics, University of Florida, Gainesville, FL 32611

\noindent
alexb\@math.ufl.edu

\Refs\tenpoint

\ref\no1\by K. Alladi\paper Partition identities involving gaps and 
weights\jour Trans. Amer. Math. Soc.\vol349\yr 1997\pages 5001-5019\endref

\ref\no2\by K. Alladi\paper A combinatorial correspondence related to 
G\"ollnitz's big partition theorem and applications\jour Trans. Amer. Math. 
Soc.\vol349\yr1997\pages2721-2735\endref

\ref\no3\by K. Alladi\paper On a partition theorem of G\"ollnitz and 
quartic transformations, ({\rm with an appendix by B. Gordon)}\jour J. Num. 
Th.\vol69\yr1998\pages153-180\endref

\ref\no4\by K. Alladi and G.E. Andrews\paper A quartic key identity for a 
partition theorem of G\"ollnitz\jour J. Num. 
Th.\vol75\yr1999\pages220-236\endref

\ref\no5\by K. Alladi, G. E. Andrews, and B. Gordon\paper  Generalizations 
and refinements of a partition theorem of G\"ollnitz\jour J. Reine Angew. 
Math.\vol460\yr1995\pages165-188\endref

\ref\no6\by K. Alladi, G. E. Andrews, and A. Berkovich\paper A new four 
parameter $q-$series identity and its partition implications ({\rm 
in preparation})\endref

\ref\no7\by K. Alladi and A. Berkovich\paper A double bounded version of 
Schur's partition theorem\jour Combinatorica - Erd\"os memorial issue 
{\rm (to appear)}\endref

\ref\no8\by K. Alladi and A. Berkovich\paper A double bounded key identity 
for a partition theorem of G\"ollnitz\jour submitted to Proc. Gainesville 
Conf. on symbolic computation (F. Garvan and M. E. -H. Ismail Eds.)\endref

\ref\no9\by K. Alladi and B. Gordon\paper Generalizations of Schur's 
partition theorem\jour Manus. Math.\vol79\yr1993\pages113-126\endref

\ref\no10\by G. E. Andrews\paper A new generalization of Schur's second 
partition theorem\jour Acta Arithmetica\vol14\yr1968\pages429-434\endref

\ref\no11\by G. E. Andrews\paper A general partition theorem with difference 
conditions\jour Amer. J. Math.\vol91\yr1969\pages18-24\endref

\ref\no12\by G. E. Andrews\paper The use of computers in the search of 
identities of Rogers-Ramanujan type\jour in Computers in Number Theory (A.O.L. 
Atkin and B. J. Birch Eds.) Academic Press\yr1971 \pages377-387\endref

\ref\no13\by G. E. Andrews\paper The theory of partitions\jour Encyclopedia 
of Math. and its Appl, Vol. 2, Addison Wesley, Reading\yr1976\endref

\ref\no14\by A. Berkovich, B. M. McCoy, and A. Schilling\paper 
Rogers-Schur-Ramanujan type identities for the $M(P,P')$ minimal models of 
conformal field theory\jour Comm. Math. Phys., \vol191\yr1998 
\pages325-395\endref

\ref\no15\by H. G\"ollnitz\paper Partitionen mit Differenzenbedingungen
\jour J. Reine Angew. Math\vol225 \yr1967\pages154-190\endref

\ref\no16\by I. Schur\paper Zur Additiven Zahlentheorie\jour Gesammelte 
Abhandlungen, Vol. 2, Springer  \yr1973\pages43-50\endref

\endRefs
\enddocument